\documentclass[12pt]{article}
\usepackage[reqno]{amsmath}
\usepackage{amssymb, theorem, enumerate}
\usepackage{graphicx}
\usepackage{tikz}
\usetikzlibrary{decorations.markings}
\usepackage{array}
\usepackage{float}
\usepackage{wrapfig}
\usepackage{fullpage}
\usepackage{url}

\expandafter\ifx\csname pdfpagebox\endcsname\relax\else
\fi
 
\theoremstyle{change}
{\theorembodyfont{\slshape}
\newtheorem{theorem}{Theorem.}[section]
\newtheorem{lemma}[theorem]{Lemma.}
\newtheorem{corollary}[theorem]{Corollary.}}
\newtheorem{conjecture}[theorem]{Conjecture.}
\theorembodyfont{\rmfamily}

\newcommand\cref[1]{Corollary~\ref{cor:#1}}

\def\proof{\noindent{{\sl Proof. }}}
\def\sqr#1#2{{\vbox{\hrule height.#2pt
    \hbox{\vrule width.#2pt height#1pt \kern#1pt
        \vrule width.#2pt}\hrule height.#2pt}}}
\def\eqed{\sqr53}
\def\qed{%
    \ifmmode\eqno\eqed
    \else\nobreak\ \hfill\eqed\medbreak\fi}

\newcommand\al{\alpha}

\newcommand\cS{{\mathcal S}}

\newcommand\Zv{{\mathbf v}}

\newcommand\ints{{\mathbb Z}}

\DeclareMathOperator\symm{Sym}
\newcommand\sym[1]{\symm(#1)}

\pgfdeclarelayer{edgelayer}
\pgfdeclarelayer{nodelayer}
\pgfsetlayers{edgelayer,nodelayer,main}
\usetikzlibrary{decorations.markings}
\tikzstyle{new}=[circle, minimum width=7pt,inner sep=2pt, fill=white,draw=black]
\tikzstyle{new2}=[circle, minimum width=7pt,inner sep=2pt, fill=black,draw=black]
\tikzstyle{none}=[circle,fill=white,draw=black,inner sep = 2pt]
\tikzstyle{n}=[shape=rectangle,minimum width=1pt,inner sep=0pt, fill=none,draw=none]
\tikzstyle{emph}=[circle, minimum width=4pt,inner sep=0pt, fill=magenta,draw=magenta]
\tikzstyle{triangle}=[rectangle, minimum width=4pt, fill=black,draw=black]
\tikzstyle{triangle2}=[rectangle, minimum width=4pt,fill=white,draw=black]
\tikzset{directed/.style={decoration={
  markings,
  mark=at position .6 with {\arrow{stealth}}},postaction={decorate}}}
\tikzset{nodirection/.style={}}

\title{Transversal polynomials of covers of graphs}
\author{Chris Godsil\footnote{University of Waterloo, Waterloo, Canada. email: \protect\url{cgodsil@uwaterloo.ca}. C. Godsil gratefully acknowledges the support of the Natural Sciences and Engineering Council of Canada (NSERC), Grant No. RGPIN-9439.}, Krystal Guo\footnote{University of Amsterdam, Korteweg-de Vries Institute, Amsterdam, The Netherlands. 
Part of this research was done when K. Guo was a post-doctoral fellow at University of Waterloo and at Universit\'{e} libre de Bruxelles, Brussels, Belgium; K. Guo gratefully acknowledges the support from ERC grant FOREFRONT
(grant agreement no. 615640) funded by the European Research Council under the EU’s
7th Framework Programme (FP7/2007-2013).  email: \protect\url{k.guo@uva.nl}}, Gordon Royle\footnote{University of Western Australia. email:\protect\url{gordon.royle@uwa.edu.au}}}
\date{December 15, 2022}

\begin{document}
\maketitle

\begin{abstract}

We explore the interplay between algebraic combinatorics and
 algorithmic problems in graph theory by defining a polynomial with
connections to correspondence colouring (also known as DP-colouring), a recent generalization of
list-colouring, and the Unique Games Conjecture. Like the chromatic
polynomial of a graph, we are able to evaluate this polynomial at a point, despite the complexity of computing this polynomial.

	We construct a cover of a graph $X$ by blowing up each vertex to a set of $r$ vertices
	and joining each pair of sets corresponding to adjacent vertices by a matching with $r$ edges.
	To each cover $Y$ of $X$ we associate a polynomial $\xi(Y,t)$, called the transversal polynomial.
    The coefficient $t^k$ of $\xi(Y,t)$ is the number of $k$-edge induced subgraphs of $Y$ whose
    vertex set is a transversal of the set system given by the blown-up vertices.
    We show that $\xi(Y,t)$ satisfies a contraction-deletion formula, and that if $n=|V_X|$ and the cover has index $r$, 
	then $\xi(Y,-(r-1)) \equiv 0 \mod r^n$. 
\end{abstract}

\section{Introduction}\label{sec:intro}

Covers of graphs occur under many guises in numerous parts of graph theory. In this paper we are concerned
with some counting problems and a particular generating function associated with covers of a graph. First
we introduce the necessary terminology.

An \textsl{arc} in an undirected graph is an ordered pair of adjacent vertices (so each edge provides 
two arcs). An \textsl{arc function of index $r$} on a graph $X$ is a function, say $\al$, from the arcs
of $X$ to the symmetric group $\sym r$, such that if $(a,b)$ is an arc in $X$, then
\[
	\al(a,b)\al(b,a) = 1.
\]
Given a graph $X = (V,E)$ and an arc function $\al$, the cover $X^\al$ is the graph with vertex set
\[
	V \times \{1,\ldots,r\}
\]
where $(x,i)$ is adjacent to $(y,j)$ if $(x,y)$ is an arc and the permutation $\al(x,y)$
maps $i$ to $j$. The sets 
\[
	\{(x,i): i=1,\ldots,r\}
\]
are called the \textsl{fibres} of the cover, they are cocliques in $X^\al$.
We may refer to a cover of index $r$ as an \textsl{$r$-fold cover}.

Colloquially speaking, we construct $X^\al$ by blowing up each vertex of $X$ into a coclique
of size $r$, and then joining two cocliques corresponding to adjacent vertices of
$X$ by a matching of size $r$. The 3-cube provides a simple example; the pairs of vertices
at distance three are the fibres of a cover of $K_4$ with index two, as shown in Figure \ref{fig:cube}. 
\begin{figure}[htb]
\centering
\begin{tikzpicture}
	\begin{pgfonlayer}{nodelayer}
		\node [style=new] (0) at (-5, 1) {};
		\node [style=new2] (1) at (-3, 1) {};
		\node [style=triangle] (2) at (-5, -1) {};
		\node [style=triangle2] (3) at (-3, -1) {};
		\node [style=n] (4) at (-4.25, -0.5) {$\epsilon$};
		\node [style=n] (5) at (-4, 0.5) {$\epsilon$};
		\node [style=n] (6) at (-2.75, 0) {$\epsilon$};
		\node [style=n] (7) at (-4, 1.25) {$(1,2)$};
		\node [style=n] (8) at (-5.25, 0) {$\epsilon$};
		\node [style=n] (9) at (-4, -1.25) {$(1,2)$};
		\node [style=new] (10) at (-1, 1.5) {\tiny $2$};
		\node [style=new] (11) at (-1.5, 1) {\tiny $1$};
		\node [style=new2] (12) at (1, 1.5) {\textcolor{white}{\tiny $2$}};
		\node [style=new2] (13) at (1.5, 1) {\textcolor{white}{\tiny $1$}};
		\node [style=triangle] (14) at (-1.5, -1) {\textcolor{white}{\tiny $2$}};
		\node [style=triangle] (15) at (-1, -1.5) {\textcolor{white}{\tiny $1$}};
		\node [style=triangle2] (16) at (1, -1.5) {\tiny $1$};
		\node [style=triangle2] (17) at (1.5, -1) {\tiny $2$};
		\node [style=triangle2] (18) at (4.75, -0.5) {};
		\node [style=triangle] (19) at (3.75, 0.5) {};
		\node [style=triangle2] (20) at (2.75, 1.5) {};
		\node [style=new2] (21) at (3.75, -0.5) {};
		\node [style=new2] (22) at (5.75, 1.5) {};
		\node [style=new] (23) at (4.75, 0.5) {};
		\node [style=triangle] (24) at (5.75, -1.5) {};
		\node [style=new] (25) at (2.75, -1.5) {};
	\end{pgfonlayer}
	\begin{pgfonlayer}{edgelayer}
		\draw [style=nodirection] (0) to (1);
		\draw [style=nodirection] (1) to (3);
		\draw [style=nodirection] (3) to (2);
		\draw [style=nodirection] (2) to (0);
		\draw [style=nodirection] (0) to (3);
		\draw [style=nodirection] (2) to (1);
		\draw [style=nodirection] (10) to (13);
		\draw [style=nodirection] (12) to (11);
		\draw [style=nodirection] (11) to (15);
		\draw [style=nodirection] (14) to (10);
		\draw [style=nodirection] (14) to (16);
		\draw [style=nodirection] (15) to (17);
		\draw [style=nodirection] (10) to (17);
		\draw [style=nodirection] (16) to (11);
		\draw [style=nodirection] (14) to (12);
		\draw [style=nodirection] (15) to (13);
		\draw [style=nodirection] (17) to (12);
		\draw [style=nodirection] (13) to (16);
		\draw [style=nodirection] (23) to (22);
		\draw [style=nodirection] (21) to (25);
		\draw [style=nodirection] (25) to (24);
		\draw [style=nodirection] (19) to (23);
		\draw [style=nodirection] (19) to (20);
		\draw [style=nodirection] (24) to (18);
		\draw [style=nodirection] (23) to (18);
		\draw [style=nodirection] (20) to (25);
		\draw [style=nodirection] (19) to (21);
		\draw [style=nodirection] (24) to (22);
		\draw [style=nodirection] (18) to (21);
		\draw [style=nodirection] (22) to (20);
	\end{pgfonlayer}
\end{tikzpicture}
\caption{An example of the cube as a $2$-fold cover of $K_4$. From left to right, we have $K_4$ with an arc-labelling function ($\epsilon$ denotes the identity permutation), the cover graph obtained, and the cover graph redrawn in the more recognizable depiction of the cube.\label{fig:cube}}
\end{figure}
A \textsl{transversal} of the set of fibres is a set of $|V|$ vertices, one from each fibre, and a subgraph of $X^\alpha$ induced by
a transversal is called a \textsl{transversal subgraph}.

Given a cover $Y$ of $X$, the object of interest to us is the generating function
\[
	\xi(Y,t) = \sum_{H} t^{|E_H|},
\]
where the sum runs over all the transversal subgraphs $H$ of $Y=X^{\alpha}$.
We call $\xi(Y,t)$ the \textsl{transversal polynomial} of the cover. If $X$ has $n$ vertices, then there are
$r^n$ transversal subgraphs and so $\xi(X^\alpha,1) = r^n$, and it is clear that $\xi(Y,t)$ is a polynomial of degree at most $|E_X|$.

There are two main results in this paper. The first is a contraction-deletion recurrence formula for $\xi(X^\alpha,t)$ 
which relies on particular variants of the usual contraction and deletion operations for graphs.  Denoting these by
by $X^{\alpha}/e$ and $X^{\alpha}\setminus e$ respectively, we get the following result:

\begin{theorem} 
	For any edge $e$ of $X$, 
	\[
		\xi(X^{\alpha},t) = (t-1) \xi(X^{\alpha}/e,t) + \xi(X^{\alpha}\setminus e, t).\qed
	\]
\end{theorem}

Using our contraction-deletion formula, we  establish a surprising congruence condition on $\xi(Y,-(r-1))$ when $Y$ is an $r$-fold cover.

\begin{theorem} 
	For any $r$-fold cover $Y$ of a graph $X$ with $n$ vertices, 
		\[\xi(Y,-(r-1)) \equiv 0 \mod r^n.\qed\] 
\end{theorem}

In Section~\ref{sec:cor-uniq}, we note some connections between the transversal polynomial and
correspondence colouring, and between the transversal polynomial and the unique games conjecture. In Section~\ref{sec:potts}, we give a connection between the transversal polynomial and the Potts model of statistical physics.

\section{Multigraphs and Covers}\label{sec:mult-covers}

Our definition of cover graph in Section \ref{sec:intro} was chosen for simplicity of exposition, but because our contraction operation may produce multiple or loops, we need more sophisticated definitions to deal with these.

First, a \textsl{looped multigraph} consists of a vertex set $V$, an edge set $E$ and a function $E$ whose values are either elements of $V$ or unordered pairs of vertices from $V$. 
If $\rho(e)$ is a single element of $V$, then $e$ is a \textsl{loop}, and otherwise $\rho(e)$ is the set of endpoints of $e$. It may occur that $\rho(e) = \rho(f)$, in which case we say that $e$ and $f$ are \textsl{parallel}. For a multigraph $X$, we will write $V_X$, $E_X$ and $\rho_X$ when the context is unclear.

We define the \textsl{valency} of a vertex to be the number of non-loop edges on it plus twice
the number of loops. We remark that this is consistent with the definition of valency in the context  of graphs embedded on surfaces.

A \textsl{homomorphism} from a multigraph $Y$ to a multigraph $X$ is a map
$h$ with domain $V_Y \cup E_Y$ that maps $V_Y$ to $V_X$ and $E_Y$ to $E_X$ such that, for any edge $e$, we have that $h(\rho_Y(e)) = \rho_X(h(e))$.  In other words, the image of the set of endpoints (in $Y$) of an edge $e$ is the set of endpoints of the image of $e$ (in $X$). As a loop has a single endpoint, loops of $Y$ are mapped to loops of $X$.

A homomorphism $h$ is a \textsl{local bijection} if for each vertex $u$, the number
of non-loop edges on $u$ and $h(u)$ are equal and the number of loops on $u$ and $h(u)$ are equal.
If $X$ and $Y$ are multigraphs and there is a locally bijective
homomorphism $h$ from $Y$ to $X$, we say that $(Y,h)$ is a \textsl{cover} of $X$. In this case
we will call $h$ a \textsl{covering map}.
If $(Y,h)$ is a cover of $X$ and $v$ is a vertex in $X$, the \textsl{fibre} of $h$ at $v$
is the set
\[
	h^{-1}(v) = \{ u\in V_Y,\ h(u)=v\}.
\]
If there is no loop on $v$, its fibre is a coclique in $Y$. If there is exactly one loop on $v$,
the fibre at $v$ is a disjoint union of cycles (and loops are cycles); in other words it is
a 2-regular graph. If there are exactly $\ell$ loops on $v$, its fibre is a $2\ell$-regular subgraph
along with a 2-factorization. If $e\in E_X$ incident to vertices $u\neq v$, its preimage relative to $h$ is a perfect matching between the fibres of $h$ and $u$ and at $v$. We will state the following lemma without proof; a proof for simple graphs can be found in Godsil \& Royle \cite{GR}.

\begin{lemma}
	If $h:Y\to X$ is a covering map and $X$ is connected, then all fibres of $h$ have the same size. \qed
\end{lemma}

We present a construction of covers of index $r$ of a looped multigraph. If $X$ is a looped multigraph, then its \textsl{arcs} are the loops of $X$ and, for each non-loop edge $e$, two oppositely-oriented directed edges between the endpoints of $e$.

An \textsl{arc function} on $X$ is a function $\al$ from the arcs
of $X$ into a group $G$, such that if $e$ and $f$ are the two 
orientations of the same edge, then
\[
	\al(e)\al(f)=1.
\]
Given a graph $X$ and an arc function $\al$ taking values in the symmetric group $\sym{r}$, we will construct the \textsl{cover of $X$ with respect to $\alpha$}, which we will denote $X^{\al}$. The cover $Y = X^{\al}$ is a multigraph with vertex set
\[
	V_X \times \{1,\ldots,r\}.
\]

For each arc $e$ from vertex $u$ to vertex $v$, we add to $Y$ the edges $e_1,\ldots, e_r$ such that 
\begin{equation}\label{matchingedge}
	\rho_{Y}(e_i) = \{(u,i),(v,i^{\al(e)})\}, \quad i=1,\ldots,r,
\end{equation}
and for each loop $f$ on the vertex $u$, we add to $Y$ the edges $f_1,\ldots, f_r$ where 
\begin{equation}\label{matchingloop}
	\rho_{Y}(f_i) = \{ (u,i),(u,i^{\al(f)}) \},\quad i=1,\ldots,r.
\end{equation}
If $f$ is a loop of $X$, then replacing $\alpha(f)$ with $\alpha(f)^{-1}$ does not alter the cover.

\begin{figure}[h]
\centering
\begin{tikzpicture}[scale=0.9]
	\begin{pgfonlayer}{nodelayer}
		\node [style=new] (0) at (-6.75, 2.25) {};
		\node [style=new] (1) at (-8.75, -0.75) {};
		\node [style=new] (2) at (-4.75, -0.75) {};
		\node [style=n] (3) at (-5.35, 1.05) {$\epsilon$};
		\node [style=n] (4) at (-5.75, 0.5) {$\epsilon$};
		\node [style=n] (5) at (-6.75, -1.4) {\footnotesize $(12)$};
		\node [style=n] (6) at (-6.75, -0.1) {\footnotesize $(12)$};
		\node [style=n] (7) at (-8.45, 1.25) {\footnotesize $(123)$};
		\node [style=n] (8) at (-7, 0.5) {\footnotesize $(132)$};
		\node [style=n] (8) at (-6.75, 3.5) {\footnotesize $(23)$};
		\node [style=new] (9) at (-0.75, 0.2) {};
		\node [style=new] (10) at (-1.5, -0.9) {};
		\node [style=new] (11) at (-0.1, -0.9) {};
		\node [style=n] (12) at (-0.75, 1.5) {$3$};
		\node [style=n] (13) at (-0.75, 2.25) {$2$};
		\node [style=n] (14) at (-0.75, 3) {$1$};
		\node [style=n] (15) at (-2.25, -1.5) {$3$};
		\node [style=n] (16) at (-2.75, -1.75) {$2$};
		\node [style=n] (17) at (-3.25, -2) {$1$};
		\node [style=n] (18) at (0.75, -1.5) {$3$};
		\node [style=n] (19) at (1.25, -1.75) {$2$};
		\node [style=n] (20) at (1.75, -2) {$1$};
		\node [style=none] (21) at (4, -1.75) {};
		\node [style=none] (22) at (6, 2.25) {};
		\node [style=none] (23) at (6, 3) {};
		\node [style=none] (24) at (8, -1.75) {};
		\node [style=none] (25) at (4.5, -1.5) {};
		\node [style=none] (26) at (7.5, -1.5) {};
		\node [style=none] (27) at (6, 1.5) {};
		\node [style=none] (28) at (3.5, -2) {};
		\node [style=none] (29) at (8.5, -2) {};
	\end{pgfonlayer}
	\begin{pgfonlayer}{edgelayer}
		\draw [style=directed, bend left = 15, looseness=1.00] (2) to (1);
		\draw [style=directed, bend left = 15, looseness=1.00] (1) to (0);
		\draw [style=directed, bend left = 15, looseness=1.00] (0) to (2);
		\draw [style=nodirection] (0) arc (270:-90:0.5);
		\draw [style=nodirection] (9) arc (270:-90:0.3);
		\draw [style=nodirection] (-0.75, 3.22) arc (270:-90:0.4);
		\draw [style=nodirection] (23) arc (270:-90:0.4);
		\draw [style=directed, bend left=15, looseness=1.00] (2) to (0);
		\draw [style=directed, bend left=15, looseness=1.00] (0) to (1);
		\draw [style=directed, bend left=15, looseness=1.00] (1) to (2);
		\draw [style=nodirection] (10) to (9);
		\draw [style=nodirection] (9) to (11);
		\draw [style=nodirection] (11) to (10);
		\draw [style=nodirection, bend left=17, looseness=1.00] (12) to (13);
        \draw [style=nodirection, bend right=17, looseness=1.00] (12) to (13);
		\draw [style=nodirection] (14) to (20);
		\draw [style=nodirection] (19) to (13);
		\draw [style=nodirection] (12) to (18);
		\draw [style=nodirection] (20) to (16);
		\draw [style=nodirection] (19) to (17);
		\draw [style=nodirection] (15) to (18);
		\draw [style=nodirection] (17) to (13);
		\draw [style=nodirection] (16) to (12);
		\draw [style=nodirection] (15) to (14);
		\draw [style=nodirection] (23) to (29);
		\draw [style=nodirection] (24) to (22);
		\draw [style=nodirection] (27) to (26);
		\draw [style=nodirection] (29) to (21);
		\draw [style=nodirection] (24) to (28);
		\draw [style=nodirection] (25) to (26);
		\draw [style=nodirection] (28) to (22);
		\draw [style=nodirection] (21) to (27);
		\draw [style=nodirection] (25) to (23);
		\draw [style=nodirection,  bend left=17, looseness=1.00] (22) to (27);
        \draw [style=nodirection,  bend right=17, looseness=1.00] (22) to (27);
	\end{pgfonlayer}
\end{tikzpicture}
\caption{A $3$-fold cover of $K_3$ with a loop. \label{fig:ex}}
\end{figure}

To illustrate these definitions we give a $3$-fold cover of a graph $X$, which is $K_3$ with a loop, in Figure \ref{fig:ex}. On the left, we show an arc function $\alpha$ on the arcs of $X$; the permutations are given in cyclic notation and $\epsilon$ denotes the identity permutation. In the center of Figure \ref{fig:ex}, we have the $X$ with the arc functions depicted as perfect matchings between sets of $\{1,2,3\}$ at each vertex. On the right of the figure, we have the graph resulting from the $3$-fold cover of $X$, as given by the arc function. We can easily compute directly that 
\[
\xi(X^{\alpha},t) = 4t^{3} + 6 t^{2} + 12 t + 5.
\]

We note that one place where covers arise in graph theory is in the study of embeddings
of graphs in surfaces. (See, e.g., Gross and Tucker \cite{GroTuc87}.) Where we use arc functions
and covers, they speak of voltage assignments and voltage graphs, and the class of covers
considered is somewhat less general than ours (for the topologically initiated, they
work with \textsl{regular covers}).

\section{Contraction and Deletion}

Many classic graph polynomial invariants, in particular the Tutte polynomial and its specializations such as the chromatic and flow polynomials satisfy some sort of contraction-deletion recurrence. In addition to providing a tool for the theoretical analysis of these polynomials, a contraction-deletion expression often forms the foundation of the only practical general purpose algorithms for computing these invariants.

Suppose that $X$ is a graph and that $e \in E_X$ is not a loop, so that $\rho(e) = \{u,v\}$. The graph $X/e$, obtained by \emph{contracting} $e$ is the graph with vertex set $V_X \setminus \{v\}$ and edge set $E_X \setminus \{e\}$ and where 
\[
\rho'(f) = \begin{cases}
\rho(f), & v \notin \rho(f);\\
\rho(f) \setminus \{v\} \cup \{u\}, & v \in \rho(f).
\end{cases}
\]
In other words, $e$ and $v$ are deleted, and any edge incident with $v$ is now incident with $u$ instead (with its other endvertex, if any, unchanged). Mentally we imagine $e$ shrinking towards $u$, dragging $v$ and any edges incident to $v$ with it, until $u$ subsumes $v$.

If $X^\alpha$ is an $r$-fold cover of $X$ and $e \in E_X$ is not a loop, so that $\rho(e) = \{u,v\}$ with $v \ne u$, then $e$ determines a matching between $h^{-1}(u)$ and $h^{-1}(v)$. Let $X^\alpha/e$ denote the graph that results when every edge of this matching is (simultaneously) contracted. The contracted edges vanish, but all other edges are dragged by their endpoints, forming loops as necessary. (This can be described formally in terms of the incidence function $\rho$, but it is analogous to the definition of $X/e$.)

It is easy to see that $X^\alpha/e$ is an $r$-fold cover of $X/e$, and the purpose of the next result is to show how to define an arc function, say $\alpha_{/e}$, so that $(X^\alpha)/e = (X/e)^{\alpha_{/e}}$. If a non-loop edge $f \in E_X$ is not incident with $v$, then it determines the same matching in $X^\alpha/e$ as in $X/e$, and so $\alpha_{/e}(a) = \alpha_{/e}(a)$ for $a \in \{f',f''\}$. Similarly if $f$ is a loop not incident with $f$, then the $2$-regular subgraph it determines within a fibre is the same in $X^\alpha$ and $X^\alpha/e$ and so again $\alpha_{/e}(f) = \alpha(f)$. 

So it is only the edges $f \in E_X \setminus \{e\}$ that are incident with $v$ that remain to be considered. 
There are three distinct cases to consider:
\begin{enumerate}
\setlength{\itemsep}{0pt}
\item The edge $f$ is not a loop and not parallel to $e$
\item The edge $f$ is parallel to $e$
\item The edge $f$ is a loop on $v$
\end{enumerate}

We must actually determine the value of $\alpha_{/e}$ on the two arcs associated with $f$. However it suffices to simply describe the value on the arc in one (specified) direction, as the value on the reverse arc is then determined.


\begin{figure}[htbp]
\begin{center}
\begin{tikzpicture}[xscale = 2, yscale=1]
\tikzstyle{vertex}=[circle, draw=black, fill = gray!10, inner sep = 0.75mm]
\tikzstyle{vertex2}=[circle, draw=black, fill = gray!50, inner sep = 0.75mm]

\foreach \x in {0,1,2} {
\foreach \y in {0,1,2} {
 \node [vertex] (v\x\y) at (\x,\y) {};
}
}

\node [above left] at (v01) {$(u,i)$};
\node [above] at (v12) {$(v,i^{\alpha(e)})$};
\node [above right] at (v20) {$(w,i^{\alpha(e)\alpha(f)}$)};

\draw [gray!40] (v00)--(v11)--(v21);
\draw [gray!40] (v02)--(v10)--(v22);

\draw [thick, -stealth, shorten >= 1mm] (v01)--(v12);
\draw [thick, -stealth, shorten >= 1mm] (v12)--(v20);

\node [vertex2] (u) at (0,-1){};
\node [vertex2] (v) at (1,-1) {};
\node [vertex2] (w) at (2,-1) {};
\node at (0,-1.5) {$u$};
\node at (1,-1.5) {$v$};
\node at (2,-1.5) {$w$};
\draw [-stealth, shorten >= 1mm] (u)--(v) node [fill=white, midway] {$e$};
\draw [-stealth, shorten >= 1mm] (v)--(w) node [fill=white, midway] {$f$};
\end{tikzpicture}
\begin{tikzpicture}[xscale = 1.5, yscale=1]
\tikzstyle{vertex}=[circle, draw=black, fill = gray!10, inner sep = 0.75mm]
\tikzstyle{vertex2}=[circle, draw=black, fill = gray!50, inner sep = 0.75mm]

\foreach \x in {0,2} {
\foreach \y in {0,1,2} {
 \node [vertex] (v\x\y) at (\x,\y) {};
}
}

\node [above left] at (v01) {$(u,i)$};
\node [above right] at (v20) {$(w,i^{\alpha(e)\alpha(f)}$)};

\draw [thick, -stealth, shorten >= 1mm] (v01)--(v20);

\draw [gray!40] (v00)--(v21);
\draw [gray!40] (v02)--(v22);

\node [vertex2] (u) at (0,-1){};
\node [vertex2] (w) at (2,-1) {};
\node at (0,-1.5) {$u$};
\node at (2,-1.5) {$w$};
\draw [-stealth, shorten >= 1mm] (u)--(w) node [fill=white, midway] {$f$};
\end{tikzpicture}
\end{center}
\caption{When $f$ is a non-loop not parallel to $e$}
\label{notlooporpar}
\end{figure}

Figure~\ref{notlooporpar} illustrates the first case, where $\rho(f) = \{v,w\}$ for some $w \notin \{u,v\}$. The left-hand diagram shows the situation in $X$ (the dark grey vertices) and $X^\alpha$ (the light gray vertices) where, by a slight abuse of notation, we re-use the labels $e$ and $f$ to denote \emph{specific arcs} associated with the edges $e$ and $f$. In particular $e$ is the arc directed from $u$ to $v$ and $f$ is the arc directed from $v$ to $w$. In $X^{\alpha}$ we see that $(u,i)$ is adjacent to $(v,i^{\alpha(e)})$ and $(v,i^{\alpha(e)})$ is adjacent to $(w,i^{\alpha(e)\alpha(f)})$. The right-hand diagram shows the situation in the contracted graphs $X/e$ and $X^\alpha/e$. After contraction, the edges associated with $f$ have endpoints $\{(u,i),(w,i^{\alpha(e)\alpha(f)})\}$. Therefore the permutation describing the matching associated with $f$ in $X^\alpha/e$ is given by $\alpha_{/e}(f) = \alpha(e) \alpha(f)$.

Now consider the second case, where we assume that $f$ has the same direction as $e$, namely directed from $u$ to $v$. So the matching associated with $f$ contains the edges $\{(u,i), (v,i^{\alpha(f)})\}$ for $1 \leqslant i \leqslant r$. In  $X^\alpha/e$, this matching becomes the set of edges $\{(u,i), (u,i^{\alpha(f) \alpha(e)^{-1}})\}$, which is a $2$-regular subgraph of the fibre associated with $u$. Expressed in terms of $X/e$, we see that $f$ is a loop on $u$ such that $\alpha_{/e}(f) =\alpha(f) \alpha(e)^{-1}$.

Finally we consider the case where $f$ is a loop on $v$, which in the contracted graph $X/e$ becomes a loop on $u$. In $X^\alpha$ the edges associated with $f$ are of the form $\{(v,i), (v,i^{\alpha(f)})\}$ for $1 \leqslant i \leqslant r$. The edge  
$\{(v,i), (v,i^{\alpha(f)})\}$ of $X^\alpha$ corresponds to the edge 
$\{(u, i^{\alpha(e)^{-1}}),(u, i^{\alpha(f) \alpha(e)^{-1}})\}$ in $X^\alpha/e$.  Putting $j= i^{\alpha(e)}$, this edge has endvertices $\{(u,j), (u,j^{\alpha(e) \alpha(f) \alpha(e)^{-1}})\}$ and so the $2$-regular subgraph is described by the permutation $\alpha_{/e}(f) = \alpha(e) \alpha(f) \alpha(e)^{-1}$.
We summarise the outcome of this discussion in a lemma. 

\begin{lemma}\label{lem:contraction-cover}
Let $X^\alpha$ be an $r$-fold cover of $X$ with arc-function $\alpha$. Let $e$ be an edge of $X$ that is not a loop. Then, using the conventions given above regarding the directions of arcs $e$ and $f$, the graph $X^{\alpha}/e$ is an $r$-fold cover of $X/e$ with arc function $\alpha_{/e}$ given by
\[
\alpha_{/e}(f) = 
\begin{cases}
\alpha(f), & f \text{ not incident to } v;\\
\alpha(e) \alpha(f), & f \text {not a loop and not parallel to } e;\\
\alpha(f) \alpha(e)^{-1}, & f \text{ parallel to } e;\\
\alpha(e)\alpha(f)\alpha(e)^{-1}, & f \text{ a loop on } v. 
\end{cases}
\]\qed
\end{lemma}

To illustrate Lemma \ref{lem:contraction-cover}, we give an example in Figure \ref{fig:contraction} of a contraction on a graph with two parallel edges. In this case $e$ is the blue edge with $\alpha(e) = (3,4)$ while $f$ is the black edge with $\alpha(f) = (1,2,3,4)$. As arcs $e$ and $f$ are both directed from left-to-right, after contraction, we have 
$
\alpha_{/e}(f) = \alpha(f) \alpha(e)^{-1} = (1,2,3,4)(3,4) = (1,2,4).
$

\begin{figure}[ht]
    \centering
\begin{tikzpicture}[scale=0.8]
	\begin{pgfonlayer}{nodelayer}
		\node [style=none] (0) at (-5, 0) {};
		\node [style=none] (1) at (-3, 0) {};
		\node [style=none] (2) at (0, 0) {};
		\node [style=none] (4) at (-5, 2) {4};
		\node [style=none] (5) at (-5, 3) {3};
		\node [style=none] (6) at (-5, 4) {2};
		\node [style=none] (7) at (-5, 5) {1};
		\node [style=none] (8) at (-3, 5) {1};
		\node [style=none] (9) at (-3, 4) {2};
		\node [style=none] (10) at (-3, 3) {3};
		\node [style=none] (11) at (-3, 2) {4};
		\node [style=none] (12) at (0, 2) {4};
		\node [style=none] (13) at (0, 3) {3};
		\node [style=none] (14) at (0, 4) {2};
		\node [style=none] (15) at (0, 5) {1};
		\node [style=none] (20) at (-5, -2) {};
		\node [style=none] (21) at (-3, -2) {};
		\node [style=none] (22) at (0, -2) {};
		\node [style=n] (23) at (-4, 0.8) {\tiny $(1234)$};
		\node [style=n] (24) at (-4, 0.25) {\tiny  $(1432)$};
		\node [style=n] (25) at (-4, -0.25) {\tiny  $(34)$};
		\node [style=n] (26) at (-4, -0.9) {\tiny  $(34)$};
		\node [style=n] (27) at (0.7, 0.5) {\tiny $(124)$};
		\node [style=n] (28) at (-4, -1.5) {\small $f$};
		\node [style=n] (29) at (-4, -2.75) {\small $e$};
		\node [style=n] (30) at (0.5, -1.7) {\small $f$};
		\node [style=n] (34) at (-5.5, -2.75) {$X$};
		\node [style=n] (35) at (0.75, -2.75) {$X/e$};
		\node [style=n] (36) at (-7.5, 0) {arc function $\alpha$};
		\node [style=n] (31) at (3, 0) {arc function $\alpha_{/e}$};
		\node [style=n] (32) at (-7, 3.5) {$X^{\alpha}$};
		\node [style=n] (33) at (3.2, 3.5) {$(X/e)^{\alpha_{/e}} \cong X^{\alpha}/e$};
	\end{pgfonlayer}
	\begin{pgfonlayer}{edgelayer}
		\draw [bend left=45] (20.center) to (21.center);
		\draw [bend right=45,blue] (20.center) to (21.center);
		\draw [style=directed, bend left=60, looseness=1.25] (0.center) to (1.center);
		\draw [style=directed,bend left=300] (1.center) to (0.center);
		\draw [style=directed,blue, bend right=60, looseness=1.25] (0.center) to (1.center);
		\draw [style=directed,blue, bend left=45, looseness=1.25] (1.center) to (0.center);
		\draw [style=nodirection] (22) arc (270:-90:0.3);
        \draw [style=nodirection] (2) arc (270:-90:0.3);
		\draw [blue] (7.center) to (8.center);
		\draw [blue] (6.center) to (9.center);
		\draw [blue, bend left=15] (5.center) to (11.center);
		\draw [blue] (10.center) to (4.center);
		\draw (7.center) to (9.center);
		\draw (6.center) to (10.center);
		\draw [bend right=15] (5.center) to (11.center);
		\draw (8.center) to (4.center);
		\draw (15.center) to (14.center);
		\draw [bend left=75, looseness=1.25] (15.center) to (12.center);
		\draw [bend left=60, looseness=1.25] (14.center) to (12.center);
		\draw [style=nodirection] (13) arc (0:360:0.4);
	\end{pgfonlayer}
\end{tikzpicture}

    \caption{An example of contracting an edge where there is a parallel edge. We contract edge $e$ in $X$, shown in blue, to obtain $X/e$. }
    \label{fig:contraction}
\end{figure}

Each edge $e$ of $X$ is associated with a set of edges of $X^\alpha$, as given by the expression \eqref{matchingedge} if $e$ is not a loop and \eqref{matchingloop} otherwise. Deleting $e$ from $X$ corresponds to deleting the associated set of edges from $X^\alpha$. This graph, which we denote $X^{\alpha} \setminus e$, is an $r$-fold cover of $X \setminus e$ with arc function $\alpha_{\setminus e}$ which is the restriction of $\alpha$ to $E \setminus e$. With contraction and deletion now well-defined, we now prove our contraction and deletion formula for the transversal polynomial. 

Each edge $e$ of $X$ is associated with a set of edges of $X^\alpha$, as given by the expression \eqref{matchingedge} if $e$ is not a loop and \eqref{matchingloop} otherwise. Therefore deleting $e$ corresponds to deleting the associated set of edges in $X^\alpha$. This graph, which we denote $X^{\alpha} \setminus e$, is an $r$-fold cover of $X \setminus e$ with arc function $\alpha_{\setminus e}$ which is the restriction of $\alpha$ to $E \setminus e$.

Of course, the main purpose of defining these graphs is to use them to find a recurrence for the transversal polynomial.

\begin{theorem} 
	\label{thm:del-con} Let $e$ be an edge of $X$ that is not a loop. 
Then
\[
\xi(X^{\alpha},t) = (t-1) \xi(X^{\alpha}/e,t) + \xi(X^{\alpha}\setminus e, t).
\]
\end{theorem}

\proof In $X^{\alpha}$, the edge $e=uv$ corresponds to the matching $\left\{ e_1, \ldots, e_r \right\}$ where $e_i$ is incident with $(u,i)$ and $(v, i^{\alpha(e)})$. We consider a transversal subgraph $S$ of $X^{\alpha}$. If $(u,i)$ and $(v, i^{\alpha(e)})$ are both vertices of $S$, for some $i$, then $S$ contains edge $e_i$ and no other edge from $\{ e_1, \ldots, e_r \}$.  Otherwise, $S$ contains none of 
$\{ e_1, \ldots, e_r\}$.

Let $\cS_0$ be the set of transversal subgraphs of $X^{\alpha}$ which contain no edge in $ \left\{ e_1, \ldots, e_r \right\}$ and  $\cS_1$ be the set of transversal subgraphs of $X^{\alpha}$ which contain exactly one edge in $ \left\{ e_1, \ldots, e_r \right\}$. Note $\cS_0$ and $\cS_1$ partition the set of all transversal subgraphs of $X^{\alpha}$ and so
\[
\xi(X^{\alpha},t) = \sum_{H \in \cS_0} t^{|E_H|} + \sum_{H \in \cS_1} t^{|E_H|}.
\]

We will show that $\cS_1$ is in one-to-one correspondence with set of transversal subgraphs of $X^{\alpha}/e$. We will then show that the set of transversal subgraphs of $X^{\alpha} \setminus e$ is in one-to-one correspondence with the disjoint union of  $\cS_0$ with the set of transversal subgraphs of $X^{\alpha}/e$. The formula will then follow.

To see that $\cS_1$ is in one-to-one correspondence with the set of transversal subgraphs of $X^{\alpha}/e$, consider $S \in \cS_1$ and let $e_i$ be the edge of $\left\{ e_1, \ldots, e_r \right\}$ in $S$. Since we contract along the perfect matching corresponding to $e$, in $X^{\alpha}/e$, the ends of $e_i$ are contracted together to form a new vertex, say $w_i$. If we replace the ends of $e_i$ with $w_i$ in $V_S$, we obtain a vertex of a transversal subgraph of $X^{\alpha}/e$, say $\psi(S)$. It is easy to see that this process can be reversed and that if $\psi(S) = \psi(S')$ then $S = S'$. We observe that $\psi(S)$ has exactly one less edge than $S$ and thus we have
\[ 
	\sum_{H \in \cS_1} t^{|E_H|} = t \, \xi(X^{\alpha}/e,t).
\]

For the next part, consider a transversal subgraph $S$ of $X^{\alpha} \setminus e$. Let $i$ be such that $(u,i)$ is a vertex of $S$ and let $j$ be such that $(v,j)$ is a vertex of $S$. If $j \neq i^{\alpha(e)}$, then the vertices of $S$ induce a transversal subgraph of $X^{\alpha}$ which does not contain any of $\left\{ e_1, \ldots, e_r \right\}$ and is thus an element of $\cS_0$. On the other hand, if  $j = i^{\alpha(e)}$ then  the vertices of $S$ induce  a transversal subgraph of $X^{\alpha}$ which contains $e_i$ and is an element of $\cS_1$. In first case, the transversal induced by the vertices of $S$ in $X^{\alpha}$ contains the same number of edges as $S$. In the second case, $S$ has the same number of edges as $\psi(S)$, with $\psi$ as given in the previous paragraph. We obtain that
\[
\xi(X^{\alpha}\setminus e, t) = \sum_{H \in \cS_0} t^{|E_H|} + \xi(X^{\alpha}/e,t).
\]
The result now follows.\qed

This contraction-deletion formula will be very useful. For example, we show in the next corollary that every $r$-fold cover of a tree on $n$ vertices has the same transversal polynomial. 

\begin{corollary} If $X$ is a tree on $n$ vertices, then any $r$-fold cover $X^{\alpha}$ of $X$ has transversal polynomial
\[\xi(X^{\alpha},t) = \sum_{j=0}^{n-1} r (r-1)^{n-1-j} \binom{n-1}{j} t^j.
\]
\end{corollary}

\proof First, we observe that since a tree has no cycles, any $r$-fold cover of $X$ is isomorphic to the disjoint union of $r$ copies of $X$. Thus we may assume that the arc function is the identity on every arc.

Letting $f_{n,r}(t)$ denote the expression given above, we observe that 
\[
f_{n,r}(t) = (t +r-1)f_{n-1,r}(t), 
\]
which we will use in combination with our deletion-contraction formula to prove the result.

We proceed by induction on the number of vertices. In the base case, if $X$ is a tree on $1$ vertex, it is easy to see that
$\xi(X^{\alpha},t) = r = f_{1,r}(t)$. If $X$ is a tree on two vertices, we can verify that 
\[ \xi(X^{\alpha},t) = rt + r(r-1) = f_{2,r}(t). 
\]

Suppose that every $r$-fold cover of a tree on $k$ vertices, where $k<n$, has $\xi(X^{\alpha},t) = f_{k,r}(t)$. We consider a tree $X$ on $n$ vertices and let $e=uv$ be a leaf of $X$ where $u$ has degree one in $X$. We will apply the contraction-deletion formula to $e$ and obtain
\[ \xi(X^{\alpha},t) = (t-1) \xi(X^{\alpha}/e,t) + \xi(X^{\alpha}\setminus e, t) = (t-1) f_{n-1,r}(t) + \xi(X^{\alpha}\setminus e, t) .\]
Observe that since vertices in the fibre of $u$ in $X^{\alpha}\setminus e$ are not incident to any edges, we have that $\xi(X^{\alpha}\setminus e, t) = r f_{n-1,r}(t)$ and so
\[ \xi(X^{\alpha},t)= (t-1 + r)f_{n-1,r}(t)
\]
and the result follows.\qed

\section{Evaluating at $-(r-1)$ modulo $r$} 
\begin{lemma}
If $X$ is a graph with $m$ loops on a single vertex, then for any $r$-fold cover $X^\alpha$
\[\xi(X^{\alpha},-(r-1)) \equiv 0 \mod r.\]
\end{lemma}

\proof 
Let $X$ be a multigraph with one vertex and $m$ loops, with arc-function $\alpha$ taking values  in $\sym{r}$. Let $n(i)$ be the number of loops $\ell$ such that $\alpha(\ell)$ has $i$ as a fixed point. We see that
\[\xi(X^{\alpha},t) = \sum_{i=1}^{r} t^{n(i)}.\]
We want to show that
\[\sum_{i=1}^{r} (-r+1)^{n(i)}\equiv 0 \mod r.
\]
We will proceed by induction on $m$, the number of loops. If $m=0$, we have that $\xi(X^{\alpha},t) =r$. Suppose $m=1$ and let $\ell$ be the only loop of $X$. Let
$F$ be the set of fixed points of $\alpha(\ell)$. We see that
\[\xi(X^{\alpha},t) = r-|F| + \sum_{i\in F} t \]
and so
\[\xi(X^{\alpha},-r+1) = r-|F| + \sum_{i\in F} -r+1 = r-|F| -r|F| + |F| = r(1-|F|) \equiv 0 \mod r.
\]
Now suppose that $X$ has at least $2$ loops and let $\ell$ be a loop of $X$. Again we let
$F$ be the set of fixed points of $\alpha(\ell)$. Let $n(i)$ be the number of loops $\ell' \neq \ell$ such that $\alpha(\ell')$ has $i$ as a fixed point. We obtain that
\[
\xi(X^{\alpha},t) = \sum_{i\notin F} t^{n(i)} +  \sum_{i\in F} t^{n(i)+1} =\sum_{i\notin F} t^{n(i)} + t \sum_{i\in F} t^{n(i)} .
\]
Observe that
\[
\xi(X^{\alpha} \setminus \ell,t) = \sum_{i\notin F} t^{n(i)} +  \sum_{i\in F} t^{n(i)}
\]
and so
\[
\xi(X^{\alpha},t)  =\xi(X^{\alpha} \setminus \ell,t) - \sum_{i\in F} t^{n(i)} + t \sum_{i\in F} t^{n(i)} = \xi(X^{\alpha} \setminus \ell,t) + (t-1) \sum_{i\in F} t^{n(i)}.
\]
We see that
\[
\xi(X^{\alpha},-r+1)=  \xi(X^{\alpha} \setminus \ell,-r+1) + (-r) \sum_{i\in F} t^{n(i)} \equiv 0 \mod r
\]
by the induction hypothesis, and the statement follows.\qed

\begin{theorem}
For any $r$-fold cover of a graph $X$ with $n$ vertices, \[\xi(X^{\alpha},-(r-1)) \equiv 0 \mod r^n.\]
\end{theorem}

\proof We proceed by (double) induction on the number of edges and the number of vertices. The inductive step follows from the contraction-deletion formula. The base case is one vertex with loops, which was proven in the previous lemma.\qed

\section{Two-fold covers}

For $r=2$ we can give a more precise statement, according to whether the underlying graph $X$ is eulerian or otherwise. Here a loop contributes $2$ to the degree of its vertex.

The following result relies on the observation that in a $2$-fold cover, any two fibres are joined by a union (possibly empty) of perfect matchings. In particular, the edges come in pairs, so that each edge $u_i v_j$ is paired with $u_{3-i} v_{3-j}$, which is the other edge in the perfect matching.

\begin{theorem}
If $X$ is an $n$-vertex graph and $\alpha$ an arc-labelling function of index $2$, then
\[\xi(X^{\alpha},-1) = \begin{cases} \pm 2^n, &\text{if }X \text{ is eulerian; and }\\
0, &\text{otherwise.}\end{cases}\]
\end{theorem}

\proof
Let $Y = X^\alpha$ be the $2$-fold cover of $X$ determined by the arc-labelling function $\alpha$. It is convenient to first assume that $X$ has no loops, so that the fibres of $Y$ are cocliques of size $2$.  If $H$ is a transversal subgraph of $Y$, then for each vertex $u \in V_X$, either $u_1$ or $u_2$ is a vertex of $H$. If $u_1$ is a vertex of $H$, then we can form a new transversal subgraph $H'$ by replacing $u_1$ by $u_2$. Consider now the effect of this exchange on the number of edges of $H$ and $H'$.  We partition the edges of $Y$ that are incident with either $u_1$ or $u_2$ into four groups according to whether the edge is incident with $u_1$ or $u_2$ and whether its other vertex is in $H$ or not. So let $A$ be the edges incident with $u_1$ and another vertex $H$, let $B$ be the edges incident with $u_1$ and a vertex not in $H$, and let $C$ and $D$ respectively be the sets of paired edges of  $A$ and $B$.

It is then clear that
\[
E_{H'} = E_H \cup D \backslash A
\]
because $A$ is the set of edges of $H$ that were lost when $u_1$ was removed, while $D$ is the set of edges that were gained when $u_2$ was included.
If $u$ (and therefore both $u_1$ and $u_2$) has degree $d$, then
$d = |A|+|B|$
and by their definition we have $|A| = |C|$ and $|B|=|D|$ and so
\[
|D| - |A| = |D| - (d - |B|) = |D| - (d - |D|) = 2 |D| - d.
\]
If $d$ is even then this number is also even, and so replacing any vertex in a transversal subgraph with the other vertex in its fibre does not change the parity of the number of edges of the corresponding transversal subgraph. Either every transversal subgraph contributes $+1$ to $\xi(X^\alpha,-1)$ or every transversal subgraph contributes $-1$ and so the result follows.

If $X$ is not Eulerian then there is some vertex $u \in V_X$ of odd degree. Then for each transversal subgraph of $X$ containing $u_1$, there is a corresponding transversal subgraph containing $u_2$ such that the numbers of edges of the two transversal subgraphs have the opposite parity and therefore cancel out in the evaluation of $\xi(X^\alpha,-1)$.

Finally we consider the situation when $X$ has loops. Any loop on the vertex $u \in V_X$ either contributes a double-edge from $u_1$ to $u_2$ or a pair of loops, one on $u_1$ and one on $u_2$.  In the former case, these ``within-fibre'' edges can never be in any transversal subgraph and they do not alter the parity of the degree of $u$, and so they can be ignored. In the latter case, if $u_1$ is replaced with $u_2$, then the loops on $u_1$ are replaced by the loops on $u_2$. Therefore the number of loops is the same for every transversal subgraph and so again they can be ignored. Insisting that each loop contributes $2$ to the degree provides a statement that incorporates graphs with loops while ensuring they have no effect on the remainder of the argument.\qed

\section{Relation to correspondence colouring and unique label cover}
\label{sec:cor-uniq}

Transversal subgraphs in covers of graphs has links to two important problems; correspondence colouring and the Unique Label Cover problem. We will explain these interesting connections in this section.

Correspondence colouring was introduced by Dvorak and Postle in \cite{DvoPos15}, as a generalization of list-colouring, to show that every planar graph without cycles of lengths $4,5,6,7,8$ is $3$-choosable. For a graph $X=(V,E)$, they define a $k$-correspondence assignment as a function $C$, which assigns to each edge $e=uv \in E$ a partial matching $C_{e}$ between $(u,1), \ldots, (u,k)$ and $(v,1), \ldots, (v,k)$. The graph is $C$-colourable if there exists a map $\phi : V \rightarrow [k]$ such that, for each edge $e=uv$, we have that $(u,\phi(u))$ and $(v,\phi(v))$ are not adjacent in $C_e$; the map $\phi$ is a correspondence colouring. The assignment is \textsl{full} at an edge if $C$ assigns to it a perfect matching. The authors of \cite{DvoPos15} note that an $r$-fold cover (also an $r$-\textsl{lift} in their terminology) is an $r$-correspondence assignment that is full at every edge.

In our setting, a correspondence colouring of $X$ with $\alpha$ as the correspondence assignment is a transversal coclique of $X^{\alpha}$ and thus the constant term of $\xi(X^{\alpha},t)$ is the number of correspondence colourings of $X^{\alpha}$.

The \emph{Unique Games Conjecture} \cite{Kho02} is important because its truth would imply the optimality of approximation algorithms for several NP-hard problems. For example, it is shown in \cite{KhoKinMosOdo07} that the algorithm of Goemans and Williamson for approximating the maximum cut in a graph, see \cite{GoeWil95}, has an approximation ratio which is optimal up to an additive constant, assuming that the Unique Games Conjecture is true and that P is not equal to NP. We will describe the problem using the terminology of this paper.

A \textsl{Unique Games instance} consists of an $r$-fold cover $X^{\alpha}$ of a graph $X$. A \textsl{labelling} $\phi$ is a mapping from $V_X$ to $[r]$. We say that a labelling $\phi$ \textsl{satisfies} an edge $e=uv$ if $(u,\phi(u))$ and $(v,\phi(v))$ are adjacent in $X^{\alpha}$. The problem is so-named because if there exists a labelling that satisfies every edge of the graph, then choosing the correct labelling at one vertex uniquely determines the labelling at the other vertices. Thus, the existence of such a labelling can be determined in polynomial time.

In this context, a labelling corresponds to a transversal subgraph of $X^{\alpha}$. The number of edges that a labelling satisfies is equal to the number of edges in the corresponding transversal subgraph. For example, a labelling satisfying every edge of $X$ is a transversal subgraph of $X^{\alpha}$ with $|E_X|$ edges. Thus, we see that the degree of $\xi(X^{\alpha},t)$ is the maximum number of edges satisfied by a labelling. We may restate the Unique Games Conjecture in our language as follows.

\begin{conjecture}[Khot \cite{Kho02}] For all $\epsilon,\delta >0$, there exists $r= r(\epsilon,\delta)$ such that given an $r$-fold cover $X^{\alpha}$ of a graph $X$, it is NP-hard to distinguish between the following two cases:
\begin{enumerate}[(a)]
\item $\deg(\xi(X^{\alpha},t)) \geq (1-\epsilon) |E_X|$;
\item $\deg(\xi(X^{\alpha},t)) \leq \delta |E_X|$.
\end{enumerate}
\end{conjecture}

\section{Connection to the Potts model} \label{sec:potts}

The Potts model arises in  statistical physics, and has close connections with the multivariate Tutte polynomial. 
Following the notation of \cite{Sok2005}, let $q$ be a positive integer and $X = (V,E)$ be a graph. We have a variable $v_e$ for each $e \in E$ and we write $\Zv = (v_e)_{e\in E}$. The \textsl{$q$-state Potts-model partition function}
for the graph $X$ is given by
\begin{equation}
   Z_G^{\rm Potts}(q, \Zv)   \;=\;
   \sum_{ \sigma \colon\, V \to \{ 1,2,\ldots,q \} }
   \; \prod_{e = uv \in E}  \,
      \biggl[ 1 + v_e \delta(\sigma_{u}, \sigma_{v}) \biggr]
   \;.
 \label{def.ZPotts}
\end{equation}
Here the sum runs over all maps $\sigma\colon\, V \to \{ 1,2,\ldots,q \}$ and $\delta$ denotes the Kronecker delta. 

The following observation is due to Guus Regts \cite{Regts}. We can consider each $\sigma$ to be a $k$-colouring (not necessarily proper) of the graph; in the product, we get a contribution of $1$ if the edge $e$ is not monochromatic in $\sigma$, and $1+v_e$ otherwise. If we set $v_e =t-1$ for some variable $t$,  then we obtain that 
\[ 
\begin{split}
   Z_G^{\rm Potts}(q, (t-1, \ldots, t-1) )   &=
   \sum_{ \sigma \colon\, V \to \{ 1,2,\ldots,q \} }
   \; \prod_{e = uv \in E}  \,
       t^{\delta(\sigma_{u}, \sigma_{v}) } \\
   &=  \sum_{ \sigma \colon\, V \to \{ 1,2,\ldots,q \} } t^{\# \text{monochromatic edges of } \sigma} \\
      &=\xi(X^{\epsilon},t)
   \end{split}
\]
where $\epsilon$ is the arc function which takes the identity element of $\sym{q}$ on every edge of $X$. 
Connections between statistical physics and the Unique Games Conjecture have also been studied in \cite{CouDavKol2020}. 

\section{Open problems and future directions}

Observe that if no correspondence colouring exists in a $r$-fold cover of a graph $X$ with arc function $\alpha$, then
\[\xi(X^{\alpha},-(r-1)) \equiv 0 \mod r-1.\]
Since $r$ and $r-1$ are consecutive integers, they are also coprime and so we see from the main theorem that
\[\xi(X^{\alpha},-(r-1)) \equiv 0 \mod r^n(r-1).\]

A natural question one can ask is: for a fixed $k \in \ints$, what are necessary and sufficient conditions such that $\xi(X^{\alpha},-(r-1)) = k r^n$? In particular, a characterization of covers $Y$ such $\xi(Y,-(r-1)) = 0$ would be interesting. 

We can also consider a bivariate version of this polynomial as follows. Let $c(Y)$ denote the number of components of $Y$ for any graph $Y$. For an $r$-fold cover $X^{\alpha}$, we define the following polynomial in variables $t$ and $s$:
\[
\zeta(X^{\alpha},t,s) =\sum_{H} t^{|E_H|}s^{c(H)},
\]
where the sum runs over all the transversal subgraphs $H$ of $X^{\alpha}$. Note that \[\zeta(X^{\alpha},t,1) = \xi(X^{\alpha},t).\] A natural question here is if the contraction-deletion formula holds in this version. 

The transversal polynomial is an invariant in the following way: if $\alpha,\beta$ are arc functions on $X$ such that there exists an isomorphism between the cover graphs which maps fibres to fibres, then $\xi(X^{\alpha},t) = \xi(X^{\beta},t) $. The converse is not true. A Hadamard matrix of order $n$ gives rise to a distance-regular, antipodal $2$-fold cover of $K_{n,n}$, called the \textsl{Hadamard graph of order $n$}; see \cite[Section 1.8]{BCN}.  There are three $20 \times 20$  Hadamard matrices which give rise to three pairwise non-isomorphic covers of $K_{20,20}$. We computed that all three covers have the same transversal polynomial, given below:
\[ \scriptsize
\begin{split}
32\, t^{160} & (1254 + 62415 t^2 + 160740 t^4 + 2346120 t^6 + 2757660 t^8 + 
   23212224 t^{10} + 18433540 t^{12} 
   + 116064255 t^{14} + 72779310 t^{16} \\
   &+
   373611060 t^{18} + 197518224 t^{20} + 883570680 t^{22} 
   +
   414563280 t^{24} + 1656586440 t^{26} + 701770320 t^{28} +
   2570801364 t^{30} \\ 
   &+ 1004575980 t^{32}
   + 3388349610 t^{34} +
   1225406520 t^{36} + 3870297720 t^{38} + 1314000936 t^{40} +
   3870297720 t^{42} \\
   &+ 1225406520 t^{44} + 3388349610 t^{46} +
   1004575980 t^{48} + 2570801364 t^{50} + 701770320 t^{52} +
   1656586440 t^{54} \\ 
   &+ 414563280 t^{56} + 883570680 t^{58} +
   197518224 t^{60} + 373611060 t^{62} + 72779310 t^{64} + 116064255 t^{66} +
   18433540 t^{68} \\
   &+ 23212224 t^{70} + 2757660 t^{72} + 2346120 t^{74} +
   160740 t^{76}  + 62415 t^{78} + 1254 t^{80}).
   \end{split}
\]
We also note a curiosity: there are four $16 \times 16$  Hadamard matrices which give rise to pairwise non-isomorphic covers of $K_{16,16}$, which give rise to four distinct transversal polynomials. This prompts the question:
If $H_1,H_2$ are Hadamard graphs of order $n$ (considered as $2$-fold covers of $K_{n,n}$), what are necessary and/or sufficient conditions for their  transversal polynomials to be equal? 
More generally, in analogy to many spectral graph theory problems, one can ask: What properties of  covers are distinguished by the transversal polynomial?


\end{document}